\documentclass[reqno]{amsart}
\usepackage{amssymb,latexsym}

\newcommand{\bfi}{\bfseries\itshape}

\newtheorem{thm}{Theorem}[section]
\newtheorem{lemma}{Lemma}[section]
\newtheorem{cor}{Corollary}[section]
\newtheorem{prop}{Proposition}[section]

\newtheorem{rem}{Remark}[section]

\makeatletter
\@addtoreset{figure}{section}
\def\thefigure{\thesection.\@arabic\c@figure}
\def\fps@figure{h, t}
\@addtoreset{table}{bsection}
\def\thetable{\thesection.\@arabic\c@table}
\def\fps@table{h, t}
\@addtoreset{equation}{section}

\makeatother

\allowdisplaybreaks
\def\intprod{\mathbin{\hbox to 6pt{%
                 \vrule height0.4pt width5pt depth0pt
                 \kern-.4pt
                 \vrule height6pt width0.4pt depth0pt\hss}}}

\begin{document}

\title{Discrete Euler-Poincar\'{e} and Lie-Poisson Equations}

\author[J.E. Marsden]{Jerrold E. Marsden}
\address{ {CDS\\California Institute of Technology, 107-81\\
Pasadena, CA 91125}}
\email{{marsden@cds.caltech.edu}}

\author[S. Pekarsy]{Sergey Pekarsky}
\address{ {CNLS,MS-B258\\Los Alamos, NM 87545}}
\address{ {CDS\\California Institute of Technology, 107-81\\
Pasadena, CA 91125}}
\email{{sergey@cds.caltech.edu}}

\author[S. Shkoller]{Steve Shkoller}
\address{ {CNLS,MS-B258\\Los Alamos, NM 87545}}
\address{ {CDS\\California Institute of Technology, 107-81\\
Pasadena, CA 91125}}
\email{{shkoller@cds.caltech.edu}}


\date{May 1998; current version January 17, 1999}

\begin{abstract}
In this paper, discrete analogues of Euler-Poincar\'{e}
and Lie-Poisson reduction theory are developed for systems on finite
dimensional Lie groups $G$ with Lagrangians $L:TG \rightarrow
{\mathbb R}$ that are  $G$-invariant.   These discrete equations provide
``reduced'' numerical algorithms which manifestly preserve the symplectic
structure. The manifold $G \times G$ is used as an approximation of $TG$,
and a discrete Langragian ${\mathbb L}:G \times G \rightarrow {\mathbb R}$
is construced in such a way that the $G$-invariance property is preserved.
Reduction by $G$ results in new ``variational'' principle for the reduced
Lagrangian $\ell:G \rightarrow {\mathbb R}$, and provides the discrete
Euler-Poincar\'{e} (DEP) equations.  Reconstruction of these equations
recovers the discrete Euler-Lagrange equations developed in \cite{MPS,WM}
which are naturally symplectic-momentum algorithms.
Furthermore, the solution of the DEP algorithm immediately leads to a
discrete Lie-Poisson (DLP) algorithm.  It is shown that when
$G=\text{SO}(n)$, the DEP and DLP algorithms for a particular choice of
the discrete Lagrangian ${\mathbb L}$ are equivalent to the Moser-Veselov
scheme for the generalized rigid body.
\end{abstract}

\maketitle

\tableofcontents

\section{Introduction}

The goal of this paper is to develop structure preserving numerical
integrators on the reduced space of a mechanical system whose
configuration space is a Lie group $G$, and whose Lagrangian
$L:TG \rightarrow {\mathbb R}$ is either left or right invariant by
the group action.  In particular, we shall develop the discrete analogue
of Euler-Poincar\'{e} theory by following the variational approach
introduced by Marsden, Patrick, and Shkoller \cite{MPS} for the
construction of discrete Euler-Lagrange equations that naturally
preserve the symplectic structure and the momentum mappings of
the Lagrangian system.

In our setting, the results of \cite{MPS} may be described as follows.
Given a Lagrangian $L:TG \rightarrow {\mathbb R}$, form the
action $S$
on curves $g:[a,b] \rightarrow G$ defined in a chart by
$$S(g(t)) = \int^b_a L(g^i(t), \dot g^i(t)) dt.$$
Allowing for arbitrary variations $\delta g$, not constrained to vanish on
$\{a,b\}$,  a computation of the first variation of $S$ leads to
\begin{equation}\label{mps1}
dS\bigl(g(t)\bigr)\cdot \delta g(t) =\int_a^b \delta g^i\left(
\frac{\partial L}{\partial g^i} - \frac{d}{dt}
\frac{\partial L}{\partial \dot g^i}\right) dt +
\left.\frac{\partial L}{\partial\dot g^i}  \delta g^i \right|_a^b.
\end{equation}
The last term of (\ref{mps1}) is a linear pairing of
$\partial L/\partial\dot g^i$, a function of $g^i$
and $\dot g^i$, with the tangent vector $\delta g^i$. Thus,
one may consider it to be a $1$-form
$\theta_{L} = (\partial {L}/\partial\dot q^i)dq^i$
on $TG$, and the
symplectic structure  is then defined by
$$ \omega_{L} = -d \theta_{L}.$$
Applying the operator $d^2$$=$$0$ to $S$, restricted to the space of
solutions of the Euler-Lagrange equations, shows that the flow $F_t$
of the Euler-Lagrange equations conserves the symplectic form; namely,
$F_t^* \omega_{L} = \omega_{L}$.
Next, let ${\mathfrak g}$
denote the Lie algebra of $G$ and define the momentum mapping
$J_\xi:TG \rightarrow {\mathbb R}$ for each $\xi \in {\mathfrak g}$
corresponding to the tangent lift of the right (or left) action of
$G$ on itself
by $J_\xi \equiv \xi_{TG} \intprod \theta_{L}$, where
$\xi_{TG}$ is the infinitesimal generator of $\xi\in {\mathfrak g}$
on $TG$.  Then, the  variational principle together with the
infinitesimal invariance of the action restricted to the space of
solutions, immediately leads to the fact that $F_t^* J_\xi = J_\xi$.
See \cite{MPS} for details.

Hence, this variational approach can be used to obtain a
symplectic-momentum integrator by discretizing $TG$ and forming a
discrete action sum.  For every choice of discretization, a unique
discrete symplectic structure is obtained, and the algorithm given by
the discrete Euler-Lagrange equations is guaranteed to preserve this
structure as well as the momentum mappings associated with it.  Our
goal is to apply the reduction procedure in this discrete setting,
restrict the Lagrangian to the reduced space,  and derive the
algorithm which preserves the induced structure.

Our procedure results in the discrete Euler-Poincar\'{e} equation
which defines an algorithm on the reduced space that is
shown to be equivalent to the discrete Euler-Lagrange equations
in the sense of reconstruction.
This reduced algorithm is used together with the coadjoint action to
advance points in $\mathfrak{g}^\ast \cong T^*G/G$ and thus approximate
the Lie-Poisson dynamics.
In subsequent papers, we shall make the extension to the more general
setting of Lagrangian reduction of a $G$-invariant system
on $TQ$ (see, for example, Cendra, Marsden, and Ratiu \cite{CMR}),
for a general manifold $Q$,
as well as to the case of dynamical systems defined on Lie algebras.

\section{The discrete Euler-Poincar\'{e} algorithm}

In this section we develop the discrete
Euler-Poincar\'{e} reduction of a Lagrangian system on $T G$.
We approximate $TG$ by $G \times G$ and form a discrete Lagrangian
${\mathbb L}: G \times G \rightarrow {\mathbb R}$ from the original
Lagrangian $L:TG \rightarrow {\mathbb R}$ as
$$
{\mathbb L}(g_k, g_{k+1}) = L (\kappa (g_k, g_{k+1}),
\chi (g_k, g_{k+1}) ),
$$
where $\kappa$ and $\chi$ are functions of $(g_k, g_{k+1})$
which approximate the current configuration $g(t) \in G$ and the corresponding
velocity $\dot{g}(t) \in T_g G$, respectively.
We choose particular discretization schemes so that the
discrete Lagrangian ${\mathbb L}$ inherits the symmetries of the original
Lagrangian $L$:  ${\mathbb L}$ is $G$-invariant on $G\times G$
whenever $L$ is $G$-invariant on $TG$.
In particular, the induced right (left) lifted action of $G$ onto $TG$
corresponds to the diagonal right (left) action of $G$ on
$G \times G$.

Having specified the discrete Lagrangian, we form the \emph{action sum}
$$
\mathbb{S} = \sum_{k=0}^{N-1} {\mathbb L}(g_k, g_{k+1})
$$
and obtain the discrete Euler-Lagrange (DEL) equations
\begin{equation}
  \label{DEL}
  D_2 {\mathbb L}(g_{k-1}, g_k) + D_1 {\mathbb L}(g_k, g_{k+1}) = 0,
\end{equation}
as well as the discrete symplectic  form $\omega_{{\mathbb L}}$ given in
coordinates on $G \times G$ by
\begin{equation}\label{dsymp}
\omega_{\mathbb L}= \frac{\partial^2 \mathbb L}{\partial g^i_k
\partial g^j_{k+1}} dg_k^i \wedge dg_{k+1}^j
\end{equation}
by extremizing $\mathbb{S} : G^{N+1} \rightarrow \mathbb{R}$ with
arbitrary variations.
It is shown in \cite{MPS} that the flow $\mathbb{F}_t$ of the
DEL equations preserves this discrete symplectic structure.
We remark here that the original canonical symplectic form $\omega$ is
also preserved by this flow. Indeed, as the discrete Legendre
transformations define a local symplectomorphism, we obtain that
$\omega (t) = F {\mathbb L}^{-1} (\omega_{{\mathbb L}}(t)) =
 F {\mathbb L}^{-1} (\omega_{{\mathbb L}}(0)) = \omega (0)$.

The discrete reduction of a right-invariant system proceeds as follows.
The induced group action on $G \times G$ is simply right multiplication
in each component:
$$
\bar g : (g_k, g_{k+1}) \mapsto (g_k \bar g, g_{k+1} \bar g),
$$
for all $\bar g, g_k, g_{k+1} \in G.$
Then the quotient map is given by
\begin{equation}\label{projection}
\pi : G \times G \rightarrow (G \times G)/G \cong G,\
\quad (g_k, g_{k+1}) \mapsto g_k g_{k+1}^{-1}.
\end{equation}
We note that one may alternatively use $g_{k+1} g_k^{-1}$ instead of
$g_k g_{k+1}^{-1}$; our choice is consistent with other literature
(see, for example, \cite{MPS}). The projection map (\ref{projection})
defines the \emph{reduced discrete Lagrangian}
$\ell : G \rightarrow \mathbb{R}$ for any $G$-invariant ${\mathbb L}$
by $\ell \circ \pi = {\mathbb L}$, so that
$$ \ell(g_k g_{k+1}^{-1}) = {\mathbb L} (g_k, g_{k+1}), $$
and the \emph{reduced action sum} is given by
$$ s = \sum_{k=0}^{N-1} \ell (f_{k k+1}), $$
where $f_{k k+1} \equiv g_k g_{k+1}^{-1}$ denote points in the quotient space.
A reduction of the DEL equations results in the {\bf discrete
Euler-Poincar\'{e}} (DEP) equations.  We state this as the following
theorem.

\begin{thm} \label{dep_thm}
Let ${\mathbb L}$ be a right invariant Lagrangian on $G \times G$, and let
$\ell: (G \times G)/G \cong G \rightarrow {\mathbb R}$ be the restriction
of ${\mathbb L}$ to $G$ given by $\ell(g_1 g_2^{-1}) = {\mathbb L}(g_1,g_2)$.
For any integer $N\ge 3$, let
$\{ (g_k, g_{k+1}) \}_{k=0}^{N-1}$ be a sequence in
$G \times G$ and define $f_{k k+1} \equiv g_kg_{k+1}^{-1}$ to be the
corresponding sequence in $G$. Then, the following are equivalent.
\begin{itemize}

\item[(1)] The sequence $\{ (g_k, g_{k+1}) \}_{k=0}^{N-1}$ is an extremum of
the action sum $\mathbb{S} : G^{N+1} \rightarrow \mathbb{R}$ for arbitrary
variations $\delta g_k = (d/d \epsilon)|_0 g_k^\epsilon$ where for each $k$,
$\epsilon \mapsto g_k^\epsilon$ is a smooth curve in $G$ such that
$g^0_k = g_k$.

\smallskip

\item[(2)] The sequence $\{ (g_k, g_{k+1}) \}_{k=0}^{N-1}$ satisfies the
{\bf discrete Euler-Lagrange equations} (\ref{DEL}).

\smallskip
\item[(3)] The sequence $\{f_{kk+1}\}_{k=0}^{N-1}$ is an extremum of the
reduced action sum $s: G \rightarrow \mathbb{R}$ with respect to
variations $\delta f_{kk+1}$, induced by the variations $\delta g_k$,
and given by
\[
\delta f_{kk+1} = T R_{f_{kk+1}} ( \delta g_k g_k^{-1} -
\operatorname{Ad}_{f_{kk+1}} \cdot \delta g_{k+1} g_{k+1}^{-1}) .
\]

\smallskip
\item[(4)] The sequence $\{ f_{k k+1} \}_{k=0}^{N-1}$ satisfies the
{\bf discrete Euler-Poincar\'{e} equations}
\begin{equation} \label{DEP}
-\ell' (f_{k-1 k}) \operatorname{Ad}_{f_{k-1 k}} T R_{f_{k-1 k}}
+ \ell' (f_{k k+1}) T R_{f_{k k+1}} = 0
\end{equation}
for $k=1,...,N-1$,
where the operators act on variations of the form
$\vartheta_k = \delta g_k g_k^{-1}$.
\end{itemize}
\end{thm}
\begin{proof}
We begin with the proof that (1) and (2) are equivalent following
\cite{MPS} and \cite{WM}.  One computes the first variation of the
discrete action ${\mathbb S}$ with variations that vanish on the set
$k=\{0,N\}$.  Thus,
\begin{eqnarray*}\nonumber
\left.\frac{d}{d \epsilon} \right|_{\epsilon=0} {\mathbb S}(g_k^\epsilon)
&=& \left.\frac{d}{d \epsilon}\right|_{\epsilon=0} \sum_{k=0}^{N-1}
{\mathbb L}(g_k^\epsilon,g_{k+1}^\epsilon)\\
&=& \sum_{k=0}^{N-1} D_1 {\mathbb L}(g_k,g_{k+1}) \delta g_k
+ \sum_{k=0}^{N-1} D_2 {\mathbb L}(g_k,g_{k+1}) \delta g_{k+1}\\
&=& \sum_{k=1}^{N-1} D_1 {\mathbb L}(g_k,g_{k+1}) \delta g_k
+ \sum_{r=1}^{N-1} D_2 {\mathbb L}(g_{r-1},g_r) \delta g_r\\
&=& \sum_{k=1}^{N-1} \left( D_1 {\mathbb L}(g_k,g_{k+1})
+ D_2 {\mathbb L}(g_{k-1},g_k) \right) \delta g_k,
\end{eqnarray*}
where we have used the discrete analogue of integration by parts which
simply shifts the sequence $g_k \mapsto g_r$ where $r=k+1$.  Since for
each $k=1,...N-1$, the variations $\delta g_k$ are arbitrary, this
establishes the DEL algorithm.  We remark that choosing variations
which do not vanish at $k=0$ and $k=N$ defines two $1$-forms whose
exterior derivative is the unique symplectic $2$-form given in
(\ref{dsymp}).

To see that (1) is equivalent to (3),  notice that since
${\mathbb L} = \ell \circ \pi$,
\[
\left. \frac{d}{d \epsilon}\right|_0 s (f^\epsilon_{kk+1}) =
\left. \frac{d}{d \epsilon}\right|_0 {\mathbb S} (g^\epsilon_k) .
\]
Now for (3) $\Leftrightarrow$ (4),  we compute
\[
\left. \frac{d}{d \epsilon}\right|_0 \sum_{k=0}^{N-1} \ell
(g^\epsilon_k {g^\epsilon_{k+1}}^{-1})
\]
 and find that
\begin{multline}
\left.\frac{d}{d \epsilon} \right|_{\epsilon=0} s(f_{kk+1}^\epsilon)
= \sum_{k=0}^{N-1} \ell'(f_{kk+1}) \left[ \delta g_k g^{-1}_{k+1} -
g_k g^{-1}_{k+1} \delta g_{k+1} g^{-1}_{k+1} \right]\\
= \sum_{k=1}^{N-1} \ell'(f_{kk+1}) \delta g_k g_k^{-1} g_k g_{k+1}^{-1}
-\sum_{r=1}^{N-1} \ell'(f_{r-1r}) g_{r-1} g_r^{-1} \delta g_r
g_r^{-1}, \nonumber
\end{multline}
where again we have used discrete integration by parts shifting the sequence
$g_k \rightarrow g_r$ with $r=k+1$, and the fact that
$\delta g_0=\delta g_N=0$.
Defining $\vartheta_k \equiv \delta g_k g_k^{-1}$, we obtain the discrete
Euler-Poincare equations (\ref{DEP}) for all variations of this form.
\end{proof}
\begin{rem}
In the case that ${\mathbb L}$ is left invariant, the discrete
Euler-Poincar\'{e} equations take the form
\begin{equation}\label{DEPleft}
-\ell'(f_{kk-1}) TR_{f_{kk-1}}
+\ell'(f_{k+1k}) \operatorname{Ad}_{f_{k+1k}}TR_{f_{k+1k}}=0,
\end{equation}
where  $f_{k+1k} \equiv g^{-1}_{k+1}g_k$ is in the left quotient
$(G\times G)/G$, and the operators act on variations of the form
$\vartheta_k = g^{-1}_k \delta g_k.$
\end{rem}

We may associate to any $C^1$ function $F$ on $G \times G$ its Hamiltonian
vector field $X_F$ satisfying $X_F \intprod \omega_{\mathbb L} = dF$. The
symplectic structure $\omega_{\mathbb L}$ naturally defines a Poisson
structure $\{ \cdot, \cdot \}_{G\times G}$ on $G\times G$ by the relation
\begin{equation}\label{p1}
\{F,H\}_{G \times G} = \omega_{\mathbb L}(X_F,X_H).
\end{equation}
\begin{thm}
If the action of $G$ on $G\times G$ is proper, then the algorithm on
$G$ defined by the discrete  Euler-Poincar\'{e} equations (\ref{DEP})
preserves the induced Poisson structure $\{ \cdot, \cdot \}_G$ on $G$
given by
\begin{equation}\label{p2}
\{f,h\}_G \circ \pi = \{ f \circ \pi, h \circ \pi\}_{G \times G}
\end{equation}
for any $C^1$ functions $f,h:(G \times G)/G \cong G \rightarrow {\mathbb R}.$
\end{thm}
\begin{proof}
Theorem 4.1 of \cite{MPS} guarantees that the DEL algorithm preserves the
symplectic structure $\omega_{\mathbb L}$ on $G\times G$; hence, by (\ref{p1}),
the DEL algorithm preserves the Poisson structure on $G \times G$.  Since
the action of $G$ on $G \times G$ is proper, the general Poisson reduction
theorem \cite{MR}
states that the projection $\pi: G \times G \rightarrow G$ is a Poisson map.

By Theorem \ref{dep_thm}, the projection of the DEL algorithm,
\[
\pi \circ (g_{k-1},g_k) \mapsto \pi \circ (g_k, g_{k+1})
\]
 is
equivalent to the DEP algorithm on G, $f_{k-1k} \mapsto
f_{kk+1}$.  Therefore, as the Poisson structure on $G$ is induced
by $\pi$ and as $\pi$ is Poisson, we have proven the theorem.
\end{proof}

As we shall prove in the following theorem,
reconstruction of the DEP algorithm (\ref{DEP}) on $G$
reproduces the DEL algorithm on $G \times G$.

\begin{thm}\label{recon_thm}
The discrete Euler-Lagrange algorithm governed by ${\mathbb L}$ and the
discrete Euler-Poincar\'{e} algorithm governed by $\ell$ are related as
follows.  The canonical projection of a solution of DEL gives a solution
of DEP, while the reconstruction of a solution of the DEP equations
results in a solution of the DEL equations.
\end{thm}
\begin{proof}
The first assertion follows by construction.  For the second
assertion, using the  definition $f_{k k+1}=g_kg_{k+1}^{-1}$, the
DEL algorithm can be  reconstructed from DEP algorithm by
\begin{equation}
(g_{k-1}, g_k) \mapsto (g_k, g_{k+1}) =
(f^{-1}_{k-1 k} \cdot g_{k-1}, f^{-1}_{k k+1} \cdot g_k),
\label{rec_DEL}
\end{equation}
where $f_{k k+1  }$ is the solution of (\ref{DEP}).
Indeed, $f^{-1}_{k k+1} \cdot g_k$ is precisely $g_{k+1}$.
Thus, at each increment, one need only compute $f_{kk+1}^{-1} \cdot g_k$
since $g_k= f^{-1}_{k-1 k} \cdot g_{k-1}$ is already known.

Similarly one shows that in the  case of a left $G$ action,
the reconstruction of the DEP equations (\ref{DEPleft})
is given by
\begin{equation}\label{ref_DEL_left}
(g_{k-1}, g_k) \mapsto (g_k, g_{k+1}) =
(g_{k-1} \cdot f^{-1}_{k k-1}, g_k \cdot f^{-1}_{k+1 k}).
\end{equation}
\end{proof}

\begin{rem}
Let us denote by $\overline{\pi}$ the quotient map
$\overline{\pi}:TG \rightarrow TG/G \cong {\mathfrak g}$ mapping
$\dot g \in T_gG$ to $\dot g g^{-1} \in {\mathfrak g}$.  In the
limit as the time step $h \rightarrow 0$, the DEL algorithm converges 
to the flow of the EL equations.

We denote the reconstruction of the flow
of the Euler-Lagrange equations from the flow of the
Euler-Poincar\'{e} equations by
${\mathfrak R}_{EP}$. Similarly, we denote the reconstruction of
the DEL algorithm from the DEP algorithm provided by Theorem
\ref{recon_thm}  by ${\mathfrak R}_{DEP}$. The following
noncommutative diagram shows these relations.

\begin{equation}\nonumber
\begin{array}{ccc}
G\times G & \stackrel{h \rightarrow 0}{\longrightarrow} & TG \\
\Big\downarrow\vcenter{%
  \rlap{$\scriptstyle{\rm }\,\pi$}} &
        & \Big\downarrow\vcenter{%
           \rlap{$\scriptstyle{\rm }\,\overline{\pi}$}}\\
G &  & {\mathfrak g}
\end{array}
\qquad \qquad
\begin{array}{ccc}
DEL & \stackrel{h \rightarrow 0}{\longrightarrow} & EL \\
\Big\uparrow\vcenter{%
  \rlap{$\scriptstyle{\rm }\,\mathfrak R_{DEP}$}} &
        & \Big\uparrow\vcenter{%
           \rlap{$\scriptstyle{\rm }\,\mathfrak R_{EP}$}}\\
DEP &  & EP
\end{array}
\end{equation}
where $G \times G \rightarrow TG$ as $h \rightarrow 0$ in the following sense.
Locally, $G \times G = F {\mathbb L}^* (T^*G)$ and as $h \rightarrow 0$,
$F {\mathbb L} \rightarrow F {L}$ which pulls-back $T^*G$
to $TG$. Thus, the  DEP
algorithm is asymptotic to the flow of the Euler-Poincar\'{e} equations
if properly interpreted by means of reconstruction.
\end{rem}

\section{The discrete Lie-Poisson algorithm}

In addition to reconstructing the dynamics on $G\times G$, we
may use the coadjoint action to form a discrete Lie-Poisson
algorithm approximating the dynamics on $\mathfrak{g}^\ast$.
Recall that in the Lie-Poisson reduction setting, for $m \in T^*_gG$,
the momentum corresponding to the velocity vector $\dot g \in T_gG$,
we define
$$
m_c = L_g^* m \in {\mathfrak g}^*,\qquad m_s = R_g^* m \in {\mathfrak g}^*
$$
to be the {\it body} and {\it spatial} momentum vectors, respectively,
with the relation
\[
m_s = \operatorname{Ad}^\ast_{g^{-1}} m_c.
\]

For the right invariant system,
the first Euler theorem states that $(d/dt)m_c=0$ (see Theorems 4.4
of Arnold and Khesin \cite{AK}), so that the body  momentum is a
constant of the motion.  For convenience, we denote the constant
$m_c$ by $\mu_0$ and $m_s(t)$ by $\mu(t)$ so that
\begin{equation}
\label{evol_mu}
\mu(t) = \operatorname{Ad}^\ast_{g^{-1}(t)} \cdot \mu_0.
\end{equation}
Now, let
${\mathcal O}\subset {\mathfrak g}$ be a coadjoint orbit.  Then
${\mathcal O}$ is a symplectic manifold with unique Kirillov-Kostant
forms $\omega^\pm$ as the coadjoint orbit symplectic structures
(see, for example, Theorem 14.4.1 in \cite{MR}).
Lemma 14.4.2 of \cite{MR} states that for any $g \in G$,
$\operatorname{Ad}^\ast_{g^{-1}} : {\mathcal O} \rightarrow {\mathcal O}$
preserves $\omega^\pm$.  On the other hand, there are  natural
Lie-Poisson $\{ \cdot, \cdot\}^\pm$ structures on ${\mathfrak g}^*$
(coming from Lie-Poisson reduction on $T^*G$) which induce  $(\pm)$
symplectic forms on each symplectic leaf in ${\mathfrak g}^*$.
These induced symplectic structures coincide with the coadjoint orbit
symplectic structures on each coadjoint orbit (see Kostant \cite{K});
hence, the coadjoint action preserves the Lie-Poisson structures.

Using the evolution equation (\ref{evol_mu}) along with the
 sequence $\{f_{kk+1}\}$ obtained by the DEP
algorithm, we find that
$$
\mu_{k+1} = \operatorname{Ad}^\ast_{g^{-1}_{k+1}} \mu_0
= \operatorname{Ad}^\ast_{(f^{-1}_{k k+1} \cdot g_k)^{-1}}  \mu_0 =
\operatorname{Ad}^\ast_{f_{k k+1}} \cdot
\operatorname{Ad}^\ast_{g_k^{-1}} \mu_0 =
\operatorname{Ad}^\ast_{f_{k k+1}} \mu_k.
$$
Thus, we have proven the following

\begin{prop}
An algorithm, called
the {\bfi discrete Lie-Poisson} (DLP) algorithm, on ${\mathfrak
g}^*$ defined 
along the sequence $\{f_{kk+1}\}$ provided by the DEP algorithm on $G$
and given by
\begin{equation}
  \label{DLP}
  \mu_{k+1} = \operatorname{Ad}^\ast_{f_{k k+1}} \cdot \mu_k
\end{equation}
is Lie-Poisson, i.e.
it preserves the $(+)$ Lie-Poisson structure on ${\mathfrak g}^*$.
\end{prop}

\begin{rem}
The corresponding discrete Lie-Poisson equations for the left invariant
system is given by \footnote{Henceforth, we shall use the notation
$\mu \in  {\mathfrak{g}}^\ast$ for the \emph{right} invariant system
and $\Pi \in  {\mathfrak{g}}^\ast$ for the \emph{left}.}
\begin{equation}\label{DLP_left}
\Pi_{k+1} = \operatorname{Ad}^\ast_{f^{-1}_{k+1 k}} \cdot \Pi_k,
\end{equation}
where $\Pi_k := \operatorname{Ad}^\ast_{g_k} \pi_0$ and
the reduced variable $m_c(t)$ is denoted by $\Pi(t)$ and
the constant $m_s$ by $\pi_0$.
\end{rem}

 Thus, one can obtain a
Lie-Poisson integrator by solving (\ref{DEP}) for $f_{k k+1}$ and
then substituting it into (\ref{DLP}) to generate the algorithm.
This algorithm manifestly preserves the coadjoint orbits and
hence the Poisson structure on
$\mathfrak{g}^\ast$. In Section \ref{sec_MV}, we shall show that
this recovers the Moser-Veselov equations for genenaralized
rigid-body dynamics on SO$(n)$.

It is instructive to compare our discrete Lie-Poisson algorithm with
that obtained by Ge and Marsden \cite{GM} using the Lie-Poisson
Hamilton-Jacobi equations. We now state their results which were
obtained for the \emph{left} action of a group $G$ on itself.
Let $H$ be a
$G$-invariant Hamiltonian on $T^\ast G$ and let $H_L$ be the
corresponding left reduced Hamiltonian on ${\mathfrak{g}}^\ast$.
If a generating function $S : G \times G \rightarrow \mathbb{R}$
of canonical transformations is invariant, then there exists a
unique function $S_L$ such that $S_L(g^{-1} g_0)= S(g,g_0)$.

The left reduced Hamilton-Jacobi equation for the function
$S_L : G \rightarrow \mathbb{R}$ is given by
\begin{equation}
  \label{LPHJ_1}
\frac{\partial S_L}{\partial t} + H_L (- T R_g^\ast \cdot
d S_L (g)) = 0,
\end{equation}
and is called the \emph{Lie-Poisson Hamilton-Jacobi} equation.
The Lie-Poisson flow of the Hamiltonian $H_L$ is generated by its
solution $S_L$;  in particular, the flow  $t \mapsto F_t$ of $S_L$
taking initial data $\Pi_0$ to $\Pi(t)$ is Poisson for each $t$ in the
domain of definition.  Next, one defines $g \in G$ as the solution of
\begin{equation}
  \label{LPHJ_2}
  \Pi_0 = - T L_g^\ast \cdot D_g S_L
\end{equation}
and then sets
\begin{equation}
  \label{LPHJ_3}
  \Pi = \operatorname{Ad}_{g^{-1}}^\ast \Pi_0.
\end{equation}
Thus, one obtains a Lie-Poisson integrator by approximately solving
(\ref{LPHJ_1}), and then using (\ref{LPHJ_2}) and (\ref{LPHJ_3})
to generate the algorithm.

Note that (\ref{LPHJ_1}) is the analogue of the usual Hamilton-Jacobi
equation
$$
\frac{\partial S}{\partial t} +
H \left( q^i, \frac{\partial S}{\partial q^i} \right) = 0
$$
and that (\ref{LPHJ_2}) and (\ref{LPHJ_3}) are the analogues of the
corresponding canonical transformations generated by a solution $S$ which
in a local chart are given by
$$
p_{0 i} = - \frac{\partial S}{\partial q_0^i}
\qquad
p_i = \frac{\partial S}{\partial q^i}.
$$

 It is interesting to
compare the Lie-Poisson Hamilton-Jacobi equation  (\ref{LPHJ_1})
with the discrete Euler-Poincar\'{e} equation (\ref{DEP}).
Although it may be less computationally intensive to solve
(\ref{LPHJ_1}),  one must further solve the difficult equation
(\ref{LPHJ_2}) before  advancing the algorithm by
(\ref{LPHJ_3}).  Therefore,  the discrete Lie-Poisson algorithm
(\ref{DLP}) provides an
\emph{explicit} and more efficient method for constructing a
Lie-Poisson dynamics as  compared to the \emph{implicit}
equations (\ref{LPHJ_2}) and (\ref{LPHJ_3}). Namely, although
both (\ref{DLP}) and (\ref{LPHJ_3}) manifestly advance the
discrete trajectory along the coadjoint orbit, the DLP equation
provides a time evolution map $\mu_k \mapsto \mu_{k+1}$ on
${\mathfrak{g}}^\ast$
using a \emph{known} solution $f_{k k+1}$, while (\ref{LPHJ_3})
advances the initial value $\Pi_0$ along the coadjoint orbit and
requires at each time step the solution $g$ of (\ref{LPHJ_2})
that approximates the current ``position'' $g(t)$.

\section{Discretization using natural charts}
In this section,
we discretize $TG$ by $G \times G$ and use
the group exponential map at the identity,
$\operatorname{exp}_e : \mathfrak{g} \rightarrow G$,
to construct an appropriate discrete Lagrangian.

\subsection{The general theory}

For finite dimensional Lie groups $G$, $\operatorname{exp}_e$
is locally a diffeomorphism and thus provides a natural chart.
Namely, there exists an open neighborhood $U$ of $e \in G$ such that
$\exp_e^{-1}: U \rightarrow {\mathfrak u}\equiv \operatorname{exp}_e^{-1}(U)$
is a $C^\infty$ diffeomorphism (this is not in general true for
infinite dimensional groups).  Hence, the manifold structure is provided
by right translation,  so that a chart at $g \in G$ is given by
\begin{equation}
  \label{chart}
\psi_g = \operatorname{exp}_e^{-1} \circ R_{g^{-1}}.
\end{equation}
We now define the  \emph{discrete Lagrangian},
${\mathbb L} : G \times G \rightarrow \mathbb{R}$, by
\begin{equation}
\label{D_Lag}
  {\mathbb L}(g_1, g_2) = {L} \left( \psi_g^{-1}
\left[ \dfrac { \psi_g(g_1) + \psi_g(g_2)}{2}\right],
(\psi_g^{-1})_*
\left[ \dfrac{\psi_g(g_2) - \psi_g(g_1)}{h}\right] \right),
\end{equation}
where $h \in \mathbb{R}_+$ is the given time step and
$g_1, g_2 \in U_g \equiv R_g(U)$.

We shall assume that $G$ has a right invariant Riemannian metric
$\langle \cdot, \cdot \rangle$ obtained by right
translating a positive bilinear form on ${\mathfrak g}$ over the
entire group.
For $K \subset G$ a compact set, we define the Riemannian distance
function, $\operatorname{dist}:K \times K \rightarrow {\mathbb R}^+$
by
\[
\operatorname{dist}(g_1,g_2) = \int_0^1 \langle \dot \gamma
(t),
\dot \gamma(t) \rangle dt ,
\]
where
$\gamma: [0,1] \rightarrow G$ is the geodesic with $\gamma(0)=g_1$ and
$\gamma(1)=g_2$.  It is then clear that diam$(U)= \operatorname{diam}(U_g)$
for all $g\in G$, so in order for (\ref{D_Lag}) to be well defined we
require that dist$(g_1,g_2) < \operatorname{diam}(U)$.  In other words, we
require that $(g_1,g_2)$ be close to the diagonal in $G \times G$.  Our
restriction on dist$(g_1,g_2)$ in turn places a restriction on the
timestep $h$.

Next, let
$$\eta = \frac{\psi_g(g_1) + \psi_g(g_2)}{2} ,$$
with corresponding group element
$$g' = \exp (\eta) \in U.$$
We denote the algebra element approximating the velocity $g^{-1} \dot{g}$
by
$$\zeta = \frac{\psi_g(g_2) - \psi_g(g_1)}{h}.$$
Using the standard formula for the derivative of the exponential
(see, for example, Dragt and Finn \cite{DF} or Channel and Scovel
\cite{CS2}) given by
$$
 T_\eta \exp = T_e R_{g'} \cdot \operatorname{iex}
(-\operatorname{ad}_\eta ), \ \ \eta \in \mathfrak{g}, \ \
g' = \exp (\eta) \in U,
$$
where iex is the function defined by
\begin{equation}\label{iex}
\operatorname{iex}(w) = \sum_{n=0}^{\infty} \dfrac{w^n}{(n+1) !},
\end{equation}
we may evaluate the push-forward of $\psi_g^{-1}$ at $\eta$.
We obtain the following expression for the discrete Lagrangian
$$
  {\mathbb L}(g_1, g_2) = {L} \left( \psi_g^{-1} (\eta),
  T_{g'} R_g \cdot T_e R_{g'} \cdot \operatorname{iex}
  (-\operatorname{ad}_\eta ) (\zeta) \right).
$$
Setting $q \equiv \psi_g^{-1} (\eta)= R_gg'$, the last formula
is expressed as
\begin{equation}
\label{D_Lag_2}
  {\mathbb L}(g_1, g_2) = {L} \left(q, T_e R_q \cdot
\operatorname{iex} (-\operatorname{ad}_\eta ) (\zeta) \right),
\end{equation}
so that locally the Lagrangian is evaluated at the base point
$q = \psi_g^{-1} (\eta) \in U_g \subset G$, and the Lie algebra
(fiber) element $\operatorname{iex} (-\operatorname{ad}_\eta)(\zeta)$
is right translated to the tangent space at the point $q$,  $T_q G$;
as $h \rightarrow 0$, this fiber element converges to the group velocity
$\dot{g} \in T_g G$.

The following lemma establishes that the discrete  Lagrangian ${\mathbb L}$
inherits the $G$-invariance property from the original
Lagrangian $L$, so that the discrete counterpart of the Euler-Poincare
reduction is well-defined.

\begin{lemma}
The discrete Lagrangian
${\mathbb L}:G \times G \rightarrow {\mathbb R}$ is right (left)
invariant under the diagonal action of $G$ on $G \times G$, whenever
${L}:TG \rightarrow {\mathbb R}$ is right (left) invariant.
\end{lemma}
\begin{proof}
We fix the right action and consider $R_{\bar g}^* ({\mathbb L})$
for some $\bar g \in G$.
By construction, $R_{\bar g} g_1, R_{\bar g} g_2 \in  R_{\bar g}(U_g)$,
whenever
$g_1,g_2 \in U_g \equiv R_g(U)$, so that the chart is given by
$\psi_{g \bar g} = \operatorname{exp}_e^{-1} \circ R_{{(g \bar g)}^{-1}}$.

By defintion, both $\eta$ and $\zeta$ are always elements of a neighborhood
of $0 \in {\mathfrak g}$, so it is clear that they are right invariant.
Hence, using the explicit form of the chart $\psi_{g \bar g}$ together
with the right invariance of the Lagrangian ${L}$, we obtain from
(\ref{D_Lag}) and (\ref{D_Lag_2}) that
\begin{eqnarray*}
\nonumber
 {\mathbb L}(R_{\bar g}g_1, R_{\bar g}g_2) &=& {L}
 \left( \psi_{g \bar g}^{-1} \left[\dfrac { \psi_{g \bar g}(g_1 \bar g)
 + \psi_{g \bar g}(g_2 \bar g)}{2}\right],
 (\psi_{g \bar g})_* \left[\dfrac{\psi_{g \bar g}(g_2 \bar g) -
 \psi_{g \bar g}(g_1 \bar g)}{h}\right] \right) \\
&=&  {L} \left( R_{\bar g} \cdot \psi_g^{-1} (\eta),
T_q R_{\bar g} \cdot   T_{g'} R_g \cdot T_e R_{g'} \cdot \operatorname{iex}
  (-\operatorname{ad}_\eta ) (\zeta) \right) \\
&=& {L} \left( R_{\bar g} \cdot q, \ T_q R_{\bar g}
\cdot T_e R_q \cdot \operatorname{iex}  (-\operatorname{ad}_\eta )
(\zeta) \right)\\
&=& {\mathbb L}(g_1, g_2). \qquad
\end{eqnarray*}
In the case that the group action is on the left, we use
$\phi_g = \operatorname{exp}_e^{-1} \circ L_{g^{-1}}$ as the chart,
 and proceed with the same argument.
\end{proof}

\begin{cor}
Using the discretization defined by (\ref{D_Lag}),
the reduced discrete Lagrangian $\ell$ determined by the projection map
(\ref{projection}), $\ell(g_1 g_2^{-1}) = {\mathbb L} (g_1, g_2)$,
can be expressed in terms of the continuous reduced Lagrangian $l$
by
\begin{equation}
\label{red_D_Lag}
\ell (g_1 g_2^{-1}) = l (\operatorname{iex} (-\operatorname{ad}_\eta )
(\zeta) ),
\end{equation}
where $\eta = (\psi_g(g_1) + \psi_g(g_2))/2$,
$\zeta = (\psi_g(g_2) - \psi_g(g_1))/h$, and $l$ can be defined
by translation to the identity of the arguments of the
right invariant Lagrangian $L$, i.e.
$l (\xi) = L (R_{g^{-1}} g, T R_{g^{-1}} \dot{g}) =
L (e, \xi)$, where $\xi = T R_{g^{-1}} \dot{g} \in
\mathfrak{g}$.
\end{cor}

The proof of this corollary follows from expression
(\ref{D_Lag_2}), and the fact that the Lagrangian $L$ is
right invariant so that translation by $q^{-1}$ to $e$
gives (\ref{red_D_Lag}).

The expressions (\ref{D_Lag_2}) and (\ref{red_D_Lag})
for the discrete Lagrangian in general
require evaluation of the infinite series for the iex function given
by (\ref{iex}); however, a simplification occurs when
$g$ is set to either $g_k$ or $g_{k+1}$.  This is due to the fact that
when $g=g_k$ or $g=g_{k+1}$, one may easily verify that
ad$_\zeta \eta := [\zeta, \eta] =0$, and hence that
$\operatorname{iex} (-\operatorname{ad}_\eta) (\zeta)=\zeta$.

For example, with $g = g_{k+1}$, the discrete Lagrangian is simply
\begin{equation}
  \label{D_Lag_3}
 {\mathbb L}(g_k, g_{k+1})={L} \left(q, T_e R_q (\zeta) \right),
\end{equation}
where
\[
\eta = \frac{1}{2} \log (g_k g_{k+1}^{-1} ), \
q \equiv \psi_{g_{k+1}}(\eta) =
(g_k g_{k+1})^{1/2}, \
\zeta = \frac{1}{h} \log (g_k g_{k+1}^{-1})
\]
and $\log \equiv \exp^{-1}$.
Consequently, the reduced discrete Lagrangian is given by
\begin{equation}
  \label{red_D_Lag_3}
\ell (f_{k k+1}) = l ( \log (f_{k k+1})/h),
\end{equation}
where $f_{k k+1} = g_k g_{k+1}^{-1}$.

Substituting the discrete Lagrangian (\ref{red_D_Lag_3}) into
the DEP equation (\ref{DEP}), we obtain the following implicit
algorithm on the Lie albebra
\begin{equation}
  \label{chart_DEP}
l' (\xi_{k k+1}/h) \cdot \chi (\operatorname{ad}_{\xi_{k k+1}}) =
l' (\xi_{k-1 k}/h) \cdot \chi (\operatorname{ad}_{\xi_{k-1 k}})
\cdot \exp (\operatorname{ad}_{\xi_{k-1 k}}),
\end{equation}
where $\xi_{k k+1} \equiv \log f_{k k+1} \in \mathfrak{g}$
and the function $\chi$ is defined to be the inverse of the
function iex defined by (\ref{iex}), $\chi(\operatorname{ad}_{\xi})
\cdot \operatorname{iex} (-\operatorname{ad}_{\xi}) =
\operatorname{Id}_{\mathfrak{g}}$. The function $\chi$ in
(\ref{chart_DEP}) arises from taking the derivative of the log
function viewed as a map from the Lie group to its algebra.
It is interesting to compare the above algorithm with the one obtained
by Channel and Scovel \cite{CS2} using the Hamilton-Jacobi
equation.

\subsection{Generalized rigid body dynamics}

We apply our DEP algorithm to the generalized rigid body problem.
In this case, $G = \operatorname{SO} (n)$ with Lie
algebra ${\mathfrak g}= {\mathfrak so}(n)$, and
the left invariant  Lagrangian is given by the kinetic energy
\begin{equation}
  \label{RB_Lag}
  L_{RB} (g, \dot{g}) = \frac{1}{2} ( \dot{g}, \dot{g} )_g =
  \frac{1}{2} \langle \dot{g}, \mathbb{J}_g ( \dot{g}) \rangle =
\frac{1}{2} \langle g^{-1}\dot{g},\mathbb{J} (g^{-1} \dot{g}) \rangle
= \frac{1}{2} (g^{-1} \dot{g}, g^{-1} \dot{g} ).
\end{equation}
Here, $\langle \cdot ,  \cdot \rangle_g$ denotes the pairing between
$T_g\operatorname{SO}(n)$ and its dual $T_g^*\operatorname{SO}(n)$
which we associate to the metric
$(\cdot, \cdot)$ on $\operatorname{SO}(n)$ by
$$ (X_g, Y_g)_g = \langle X_g, {\mathbb J}_g Y_g \rangle_g, \ \
X_g, Y_g \in T_g\operatorname{SO}(n),$$
where ${\mathbb J}_g = (L_{g}^\ast)^{-1} \ \mathbb{J} \ (L_{g^{-1}})_\ast$
is the left translated inertia tensor, and
${\mathbb J}: {\mathfrak so}(n) \rightarrow {\mathfrak so}(n)^*$.
On $\operatorname{SO}(n)$, $(L_{g^{-1}})_\ast \cdot \dot{g}=g^{-1}\dot{g}$.

We discretize $T^\ast \operatorname{SO} (n)$ by
$\operatorname{SO} (n) \times \operatorname{SO} (n)$ and construct
the discrete Lagrangian following (\ref{D_Lag_3}) as
$$
 {{\mathbb L}}_{RB}(g_k, g_{k+1})=L_{RB}
\left(q_{k+1 k}, T_e L_{q_{k+1 k}} (\zeta_{k+1 k})\right),
$$
where $q_{k+1 k} = g_{k+1} (g^{-1}_{k+1} g_k)^{1/2}$ and
$\zeta_{k+1 k} = \frac{1}{h} \log (g_{k+1}^{-1} g_k )$.
Using  the left invariance of the metric,
we may express the discrete rigid body Lagrangian as
\begin{equation}
  \label{D_RB_Lag}
  {{\mathbb L}}_{RB}(g_k, g_{k+1})
= \frac{1}{2} (\zeta_{k+1 k},  \zeta_{k+1 k})
= \frac{1}{2} \langle \zeta_{k+1 k}, \mathbb{J} (\zeta_{k+1 k}) \rangle.
\end{equation}

The Lagrangian for the reduced system on $(\operatorname{SO}(n) \times
\operatorname{SO}(n)) /\operatorname{SO}(n) \cong \operatorname{SO}(n)$
is then given by
\begin{equation}
  \label{D_RB_Lag_red}
  {\ell}_{RB} (f_{k+1 k})={{\mathbb L}}_{RB}(g_k, g_{k+1})=\frac{1}{2 h^2}
  \langle \log f_{k+1 k}, \mathbb{J} (\log f_{k+1 k}) \rangle,
\end{equation}
where $f_{k+1 k} \equiv g_{k+1}^{-1} g_k \in \operatorname{SO} (n)$
is an element of the reduced space and $h$ is the time step.

The DEP equation (\ref{DEPleft}) has the following implicit form
\begin{equation}
  \label{DEP_RB}
\zeta_{k+1 k} = \mathbb{J}^{-1} \left( \operatorname{iex}
(- \operatorname{ad}^\ast_{h \zeta_{k+1 k}}) \cdot
\chi ( \operatorname{ad}^\ast_{h \zeta_{k k-1}}) \cdot
\operatorname{Ad}^\ast_{\exp (- h \zeta_{k k-1})}
\mathbb{J} (\zeta_{k k-1}) \right).
\end{equation}

\section{Moser-Veselov discretization of the generalized rigid body}
\label{sec_MV}

An alternative discretization approach may be taken if we
first embed our group $G$ into a linear space; for finite dimensional
matrix groups, the linear ambient space is
$\mathfrak{gl}(n)$.
Then, summation of the group elements becomes a legitimate operation
provided we project the result back onto the group $G$ by using
Lagrange multipliers.

In this section, we consider
the \emph{left} invariant generalized rigid body equations on
$\operatorname{SO}(n)$. The corresponding Lagrangian is
determined by a symmetric positive definite operator
$J : \mathfrak{so}(n) \rightarrow \mathfrak{so}(n)$, defined by
$J (\xi) = \Lambda \xi + \xi \Lambda$, where $\xi \in \mathfrak{so}(n)$
and $\Lambda$ is a diagonal matrix satisfying $\Lambda_i + \Lambda_j>0$
for all $i \ne j$. The left invariant metric on
$\operatorname{SO}(n)$ is obtained by left translating the bilinear
form at $e$ given by
$$
( \xi, \xi ) = \frac{1}{4} \operatorname{Tr} \left( \xi^T J (\xi)
\right).
$$

The operator $J$, viewed as a mapping $\mathbb{J} : \mathfrak{so}(n)
\rightarrow \mathfrak{so}(n)^\ast$, has the usual interpretation of
the inertia tensor, and the
$\Lambda_i$ correspond to the sums of certain principal moments of
inertia.

The rigid body Lagrangian is the kinetic energy of the
system
\begin{equation}
\label{grb_Lag}
L (g, \dot{g}) = \frac{1}{4} \langle  g^{-1} \dot{g},
\mathbb{J} ( g^{-1}\dot{g}) \rangle = \frac{1}{4} \langle \xi,
\mathbb{J}(\xi) \rangle,
\end{equation}
where $\xi =  g^{-1}\dot{g} \in \mathfrak{so}(n)$ and
$\langle \cdot , \cdot \rangle$
is the pairing between the Lie group and its dual; hence, the Hamiltonian
vector field of $L$ is the geodesic spray on $TG$.

Using the definition of $J$ we rewrite the Lagrangian
(\ref{grb_Lag}) in the following form:
$$
L = \frac{1}{4} \operatorname{Tr} \left( \xi^T J (\xi)
\right) = \frac{1}{2}\operatorname{Tr} \left( \xi^T \Lambda \xi \right).
$$
We now discretize the Lie algebra elements by
$\xi = g^{-1} \dot{g}$
\begin{equation}
\label{MV_disc}
\xi \approx \frac{1}{h} g_{k+1}^T (g_{k+1} - g_k),
\end{equation}
where $h$ is the time step. Substituting (\ref{MV_disc}) into the Lagrangian
$L$ (and using properties of the trace),
we obtain the following expression for the discrete Lagrangian
(modulo ${\mathbb R}$):
$$
{\mathbb L}(g_k, g_{k+1}) = - \frac{1}{h^2} \operatorname{Tr} \left(
g_k \Lambda g^T_{k+1} \right).
$$
We remark that exactly the same expression is obtained if we instead
discretize $\xi$ by $\frac{1}{h} g_k^T(g_{k+1}-g_k)$.
Notice that up to a multiplier of $-1/h^2$, this is
precisely the Lagranigian used by Moser and Veselov \cite{MoV}.

We scale the above Lagrangian and introduce matrix Lagrange multipliers
$\lambda_k$, imposing the constraint $\Phi_k(g_k) = g_k g_k^T -
\operatorname{Id} = 0$.
By decomposing $\lambda_k$ into symmetric and skew components,
we see that the skew component of $\lambda_k$ does not contribute to
the action because the constraint $\Phi_k$ is symmetric; thus, we
find that $\lambda_k = \lambda_k^T$.
The action sum then takes the form
\begin{equation}
\label{grb_action}
S = \sum_k \operatorname{Tr} \left(g_k \Lambda g^T_{k+1} \right)
- \frac{1}{2} \sum_k \operatorname{Tr} \left( \lambda_k (g_k g_k^T -
\operatorname{Id}) \right)
\end{equation}

Notice that the discrete Lagrangian ${\mathbb L}$ is left invariant and
can be reduced to a Lagrangian $\ell : G \rightarrow \mathbb{R}$
using the canonical projection $\pi : (g_k, g_{k+1}) \mapsto
f_{k+1 k} = g_{k+1}^{-1} g_k$ so that
$$ \ell (f_{k+1 k}) = \operatorname{Tr} (f_{k+1 k} \Lambda).$$
Because the constraint, ensuring that each $g_k \in G$, is $G$-invariant,
there exists a Lagrange multiplier $\bar\lambda_k$ in the conjugacy
class of $\lambda_k$, i.e., $\bar\lambda_k = g^T \lambda_k g$ for all
$g \in G$, so that $\bar \lambda_k = \bar \lambda_k^T$.
Hence, computing the discrete variation of $\operatorname{Tr} \left( \lambda_k
 \Phi_k(g_k) \right)$ with respect to $g_k$, we obtain the operator equation
$$
-\ell'(f_{kk-1}) TR_{f_{kk-1}}
+\ell'(f_{k+1k}) \operatorname{Ad}_{f_{k+1k}}TR_{f_{k+1k}} =
\bar\lambda_k ,
$$
where the operators act on the variations $\vartheta_k = g_k^T \delta g_k$.
Using the expression for the reduced Lagrangian $\ell$, the DEP
equation can then be written as
$$
f_{k+1 k}^T \Lambda + f_{kk-1} \Lambda = \bar \lambda_k .
$$
Using the fact that $\bar\lambda_k^T = \bar\lambda_k$, we
obtain the DEP algorithm on $\operatorname{SO}(n)$ as
\begin{equation}
\label{grb_DEP}
f_{k+1 k}^T \Lambda - \Lambda f_{k+1 k} =
\Lambda f_{kk-1}^T - f_{kk-1} \Lambda.
\end{equation}
This is an implicit scheme to be solved for $f_{k+1 k}$ using
the current value $f_{kk-1}$. The solution of (\ref{grb_DEP}) generates the
\emph{explicit} DLP algorithm on $\mathfrak{so}(n)^*$ given
by
\begin{equation}
\label{grb_DLP}
\Pi_{k+1} = \operatorname{Ad}_{f_{k+1 k}^{-1}} \Pi_k =
f_{k+1 k} \Pi_k f_{k+1 k}^T.
\end{equation}

Finally, reconstruction of the DEP algorithm recovers the DEL algorithm
on $G\times G$ which, according to (\ref{rec_DEL}), is given by
$$ (g_{k-1}, g_k) \mapsto (g_k, g_{k+1}) = (g_k, g_k \cdot f^{-1}_{k+1 k}).$$

\begin{thm}\label{thm5.1}
The above DEP and DLP algorithms given by (\ref{grb_DEP}) and
(\ref{grb_DLP}), respectively, are equivalent to the Moser-Veselov
equations
\begin{equation}
\label{grb_eqn}
\left\{ \begin{array}{l}
M_{k+1} \equiv \omega_{k-1} M_k \omega_{k-1}^{-1} \\
M_k = \omega_k^T \Lambda - \Lambda \omega_k, \qquad
\omega_k \in \operatorname{SO}(n),
\end{array} \right.
\end{equation}
where (using the notation of \cite{MoV})
$\omega_k = g_k^T g_{k-1} \in \operatorname{SO}(n)$ is the
discrete angular velocity, $M_k = g_{k-1}^T m_k g_{k-1} =
\omega_k^T \Lambda - \Lambda \omega_k \in \operatorname{so}(n)$ is the
discrete body angular momentum, and $m_k = m_0$
is the constant discrete spatial angular momentum.
\end{thm}
\begin{proof}
Comparing the definitions of $f_{k k-1} = g_k^T g_{k-1}$
and $\omega_k = g_k^T g_{k-1}$,
we see that $f_{k k-1} \equiv \omega_k$.
Similarly, comparing the definitions of $\Pi_k =
\operatorname{Ad}^\ast_{g_k} \pi_0$ and
\[
M_k = g_{k-1}^T m_k g_{k-1}
= g_{k-1}^T m_o g_{k-1} = \operatorname{Ad}^\ast_{g_{k-1}} m_0,
\]
 we
conclude that $\Pi_{k-1} \equiv M_k$ and $\pi_0 \equiv m_0$.
Hence,  the first equation in (\ref{grb_eqn}) is
precisely the DLP algorithm (\ref{grb_DLP}).

Substituting the second equation of (\ref{grb_eqn}) into the first
results in the following expression:
$$
\omega_{k+1}^T  \Lambda - \Lambda \omega_{k+1}
= \Lambda \omega_k^T - \omega_k  \Lambda,
$$
which is precisely the DEP equation (\ref{grb_DEP}) when the above
identifications are invoked.
\end{proof}

The Moser-Veselov algorithm (\ref{grb_eqn}) has an an
obvious geometric mechanical interpretation.
The first equation can be viewed as a discretization of
the left Lie-Poisson equation
$$M_k = g_{k-1}^T m_0 g_{k-1} = \operatorname{Ad}_{g_{k-1}}^\ast m_0,$$
rewritten in terms of the $\omega_k$ and this corresponds to the DLP
algorithm (\ref{grb_DLP}).
The second equation is a discrete version of the relation between
the angular momentum and angular velocity, as it is obtained by
substitution of (\ref{MV_disc}) into $M = J (\xi) = \Lambda \xi +
\xi \Lambda$.

The DEP algorithm (\ref{grb_DEP}) provides an equivalent alternative
to the Moser-Veselov scheme (\ref{grb_eqn}), the difference
being that the former is an algorithm on $G$ only, while the
latter is a combined algortihm on $G$ and $\mathfrak{g}^\ast$
and schematically can be represented by the mappings
$\mathfrak{g}^\ast \mapsto G \mapsto \mathfrak{g}^\ast \mapsto G;
\quad M_k \mapsto \omega_k \mapsto M_{k+1} \mapsto \omega_{k+1}$.

In proof of Theorem \ref{thm5.1}, we identified $\Pi_{k-1}$ with
$M_{k}$ in order to establish the equivalence with the Moser-Veselov
algorithm;  however, without any such identification, we exactly obtain
the algorithm  given by equation (4.1) in Lewis and Simo \cite{LS} which
we write in our notation as
\begin{eqnarray}
g_{k+1}&=& g_k f_{k+1k}^T, \nonumber\\
\Pi_{k+1}& =& f_{k+1k}\Pi_k f_{k+1k}^T,\label{ls2}\\
\Delta t \Pi_k& =& 2 \operatorname{skew}(g_k \Lambda). \nonumber
\end{eqnarray}
Equation (\ref{ls2}$_1$) corresponds to our reconstruction algorithm
(\ref{ref_DEL_left}), (\ref{ls2}$_2$) corresponds to our DLP algorithm
(\ref{DLP_left}),  and (\ref{ls2}$_3$) is our DEP algorithm
(\ref{grb_DEP}).  To see this, simply note that
$$g_{k}^T \left( \text{\cite{LS}, Eq. 4.5} \right) g_k = 
\text{Eq. \ref{grb_DEP}  (i.e. DEP}).$$

It is worthwhile to make a few remarks at this point.  Although it is
claimed in \cite{LS} that a computation of the first variation of the
action $\sum_{k} \text{Tr}(g_k \Lambda g_{k+1}^T)$ leads to the
algorithm (\ref{ls2}), we have shown that only constrained variations
of the action function (\ref{grb_action}) lead to this algorithm.
Furthermore, the algorithm (\ref{ls2}) is obtained by constraining
the iterates of the momentum to be equal; this constraint is superfluous
as the discrete Euler-Lagrange equations necessarily conserve the momentum.
Finally, if we choose $f_{k+1k} = \text{cay}(\xi_{k+1k})$ where
cay$:\text{so}(n)
 \rightarrow \text{SO}(n)$ is the Cayley transform given by
cay$(\xi) = (1+{\frac{1}{2}}\xi) (1-{\frac{1}{2}}\xi)^{-1}$ for any
$\xi \in \text{so}(n)$, then the rigid-body algorithm for $\xi_{k+1k}$
is second-order accurate, as proven in \cite{LS}.  It is not clear, however,
whether the second-order accuracy can be maintained in the absence of the
Cayley transform.

\section{A comparison of DEP/DLP algorithms with splitting methods}
For the purpose of comparison, we shall now describe the Hamiltonian
splitting methods for generating Lie-Poisson integrators on ${\mathfrak g}^*$,
the dual of the Lie algebra of a group $G$.  The basic idea behind
the construction of such an
algorithm follows from the fact that many Lie-Poisson systems are governed by
reduced Hamiltonians $h$ which can be written as a
sum $h^1 + \dots + h^N$, where each $h^i$
can be exactly integrated.  Letting $\phi^i_t$ denote the flow of the
Hamilonian system $h^i$, we see that to first order in the time-step
$\Delta t$,  the flow $\phi_t$ generated by $h$ may
be expressed as
$$ \phi_{\Delta t} = \phi^1_{\Delta t} \circ \dots \circ \phi^N_{\Delta t}.$$
As each of the maps $\phi^i_{\Delta t}$ is a Poisson map, hence symplectic
on each leaf, the composition must also preserve the Poisson structure.
Consequently,
all Casamirs are also preserved by this splitting algorithm.  Furthermore,
one may construct this splitting algorithm to any order of accuracy in
$\Delta t$.
(For example, the leapfrog method $\phi_{{\frac{1}{2}}\Delta t}
\phi^{-1}_{-{\frac{1}{2}}\Delta t}$ is a second order accurate scheme (see,
for example, \cite{MS1}).)

Whereas the DEP/DLP algorithms manifestly preserve the Poisson structure and
all of the corresponding Casamirs as well, they do much more.   First, the
reduced algorithms may be used in both the Lagrangian and Hamiltonians sides,
in that computation of the discrete Euler-Poincar\'{e} trajectory immediately
leads to the discrete Lie-Poisson trajectory on ${\mathfrak g}^*$.
More importantly, the discrete Lie-Poisson or Euler-Poincar\'{e} dynamics
may be
reconstructed to obtain symplectic-momentum integrators on $TG$, for example.
Conservation of momentum ensures
that the reconstructed discrete trajectory lies in an $n$ dimensional
submanifold
of the full $2n$ dimensional space $G \times G$, approximating $TG$.  This $n$
dimensional submanifold is the level set of the discrete momentum mapping.  For
a small enough time step $\Delta t$, $G \times G$ is locally diffeomorphic
to $TG$
through the discrete Legendre transform, and hence we ensure that our discrete
reconstruced trajectory is conserving the actual momentum.

Now recall that for right invariant systems, we have used the variable $m_s$ to
denote the solution of the Lie-Poisson equation, from which we obtain that
$m_c(t) \equiv \operatorname{Ad}^*_{g(t)} m_s(t)$ is conserved.  Using our
DEP algorithm, we may compute the discrete trajectory $\{ {(m_s)}_{kk+1}\}$,
reconstruct to find $g_{k}$, and find that
${(m_c)}_{kk+1}=\operatorname{Ad}^*_{g_k}
{(m_s)}_{kk+1}$ is conserved.  On the other hand, the splitting method does
not provide an algorithm for reconstructing the motion on $T^*G$ in such a
way as
to ensure conservation of momentum;  thus, there is no obvious way to define
the discrete analogue of $m_c$, let alone check that it is conserved.

Nevertheless, there are some computational advantages to using the
splitting method;
the fact that the splitting method leads to an explicit scheme is perhaps
the most
important of these advantages. An efficient excplicit algortihm for the
SU$(n)$ model of two dimensional hydrodynamics on a torus is constructed
in \cite{Mc}. The author presents a Poisson integrator of complexity
O$(N^3 \log N)$ which preserves $N-1$ Casimirs.

\section{Addendum: relation to other works}

It is very interesting to compare the above constructions
and algorithms to the recent results of Bobenko and Suris
\cite{BS}. In this paper they consider the theory of discrete 
time Lagrangian mechanics on Lie groups and, more specifically,
address the issue of discrete Lagrangian reduction using
left or right trivializations of the (co)tangent bundles 
of Lie groups. They adopt a somehow broader point of view
when the symmetry group of a system defined on a Lie group $G$
is a subgroup of $G$. Hence, it includes the Lie-Poisson case as
a special case. Below we shall demonstrate that the reduced
discrete equations obtained in \cite{BS} agree with our
DEP/DLP algorithms when the symmetry group is taken to be
the full group $G$.
Here we summarize their results  choosing for consistency
and simplicity the case of right trivialization and refer the
reader to \cite{BS} for details of proofs and notations.

Let the discrete Lagrangian ${\mathbb L}(g_k, g_{k+1}) :
G \times G \rightarrow {\mathbb R}$ define a a discrete system
with the corresponding DEL equations.
Consider the map 
\begin{equation}
  \label{bs_rmap}
(g_k, w_k) \in G \times G \ \mapsto 
\ (g_k, g_{k+1}) \in G \times G,
\end{equation}
where 
$$
g_{k+1}=w_k g_k \ \Leftrightarrow \ w_k = g_{k+1} g_k^{-1}.
$$
Consider also the right trivialization of the cotangetn bundle
$T^* G$:
$$
(g_k, m_k) \in G \times {\mathfrak g}^* \ \mapsto 
\ (g_k, \Pi_k) \in T^* G,
$$
where
$$
\Pi_k = R_{g_k^{-1}}^* m_k \ \Leftrightarrow \
m_k = R_{g_k}^* \pi_k.
$$

Denote the pull-back of the Lagrange function under (\ref{bs_rmap})
by
$$
{\mathbb L}^{(r)} (g_k, w_k) = {\mathbb L}(g_k, g_{k+1}).
$$
Proposition $3.5$ of \cite{BS} gives the DEL equations in these 
coordinates:
\begin{equation}
\label{bs_DEL}
\left\{ \begin{array}{l}
\text{Ad}^*_{w_k} m_{k+1} = m_k + d_g {\mathbb L}^{(r)} (g_k, w_k), \\
g_{k+1} = w_k g_k,
\end{array} \right.
\end{equation}
where
$$
m_k = d_w {\mathbb L}^{(r)} (g_{k-1}, w_{k-1}) \in {\mathfrak g}^*.
$$

Assume that for some $\zeta \in {\mathfrak g}$, ${\mathbb L}^{(r)}$ 
is invariant under the action of a subgroup 
$G^{[\zeta]} \equiv \{ h \in G | \text{Ad}_h \zeta = \zeta \} \subset G$
on $G \times G$ induced by right translations on $G$:
$$
{\mathbb L}^{(r)} (g h, w) = {\mathbb L}^{(r)} (g, w), 
\quad h \in G^{[\zeta]}.
$$
Define the reduced Lagrange function 
$\Lambda^{(r)} : {\mathfrak g}_\zeta \times G \mapsto {\mathbb R}$
as
$$
\Lambda^{(r)} (a, w) = {\mathbb L}^{(r)} (g, w), \quad
a = \text{Ad}_g \zeta \in {\mathfrak g}_\zeta
$$
here ${\mathfrak g}_\zeta$ is the adjoint orbit of $\zeta$.

Then, proposistion $3.7$ of \cite{BS} states that under the
reduction by $G^{[\zeta]}$, the reduced Euler-Lagrange equations
become
\begin{equation}
\label{bs_redEL}
\left\{ \begin{array}{l}
\text{Ad}^*_{w_k} m_{k+1} = m_k - 
\text{ad}^*_{a_k} \nabla_a \Lambda^{(r)} (a_k, w_k), \\
a_{k+1} = \text{Ad}_{w_k} a_k,
\end{array} \right.
\end{equation}
where
$$
m_k = d_w \Lambda^{(r)} (a_{k-1}, w_{k-1}) \in {\mathfrak g}^*.
$$
In (\ref{bs_redEL}) the follwoing notations are used (see \cite{BS}).
For a function $f: G \mapsto {\mathbb R}$, 
its left and right Lie derivatives,
$d f(g) : G \mapsto {\mathfrak g}^*$ and 
$d' f(g) : G \mapsto {\mathfrak g}^*$, are defined via
$$
\langle d f(g), \eta \rangle = \frac{d}{d \epsilon} 
f (\text{e}^{\epsilon \eta} g) |_{\epsilon = 0}, 
\quad \forall \ \eta \in {\mathfrak g},
$$
$$
\langle d' f(g), \eta \rangle = \frac{d}{d \epsilon} 
f (g \text{e}^{\epsilon \eta}) |_{\epsilon = 0}, 
\quad \forall \ \eta \in {\mathfrak g}.
$$
Then, the gradiaent $\nabla f : G \mapsto T^* G$ is related 
to the above derivatives via
$$
\nabla f(g) = R_{g^{-1}}^* d f(g) = L_{g^{-1}}^* d' f(g).
$$

Notice that in the Lie-Poisson case, when the symetry group is 
$G$ itself, the reduced space is simply the group $G$ represented 
by $w_k = g_{k+1} g_k^{-1}$, and equations (\ref{bs_redEL})
become
\begin{equation}
\label{bs_DEP}
\text{Ad}^*_{w_k} m_{k+1} = m_k 
\end{equation}
with 
\begin{equation}
  \label{bs_mom}
m_k = d \Lambda^{(r)} (w_{k-1}) \in {\mathfrak g}^*.
\end{equation}

Comparing the above notations with the results in our paper,
we immediately see that $w_k$ correspond to the other choice
for the quotient map (\ref{projection}) 
$\pi : (g_k, g_{k+1}) \mapsto f_{k k+1} \equiv g_k g_{k+1}^{-1}$,
i.e. $f_{k k+1} = w_k^{-1}$.
Similarly, the reduced Lagrangian $\Lambda^{(r)} (w_k)$
corresponds to $\ell (f_{k k+1})$ in our notations.
Finally, using the definitions of the Lie derivatives above we 
obtain for the angular momentum (\ref{bs_mom})
$$
m_k = d \Lambda^{(r)} (w_{k-1}) = R^*_{w_{k-1}} 
\nabla \Lambda^{(r)} (w_{k-1}) 
= R^*_{f_{k-1 k}^{-1}} \ell' (f_{k-1 k}),
$$
where we have substituted our notations.
Hence, (\ref{bs_DEP}) can be written as
$$
{Ad}^*_{f_{k k+1}^{-1}} R^*_{f_{k k+1}^{-1}} \ell' (f_{k k+1})
= R^*_{f_{k-1 k}^{-1}} \ell' (f_{k-1 k}).
$$

The last expression is precisely the DEP algorithm (\ref{DEP})
after rewriting it with the adjoints of the above operators acting
on the variation $\vartheta_k = \delta g_k g_k$ (see section $2$).
It is interestig to note that the second equation in (\ref{bs_DEL})
corresponds to our reconstruction equation (\ref{rec_DEL}).
Similar correspondence can be estableshed for the case of left
trivialization considered in \cite{BS}.

\section*{Acknowledgments}
The authors would like to thank Anthony Bloch, Peter Crouch, and
Tudor Ratiu for helpful comments.


\end{document}